\documentclass{article}
\usepackage[utf8]{inputenc}
\usepackage[margin =2.4cm]{geometry}
\usepackage{amssymb}
\usepackage{amsmath}
\usepackage{amsthm}
\usepackage{enumerate} 
\usepackage{graphicx}
\usepackage{csquotes}
\usepackage{subcaption}
\usepackage{multirow}
\usepackage{url}

\title{Selection of Criteria Using MCDM Techniques - An Application in Renewable Energy}
\author{Mahak Bhatia \\
Department of Science \\ 
St. Xavier's College \\ 
Jaipur, India \\
\texttt{drmahakbhatia@gmail.com}
 \and Aled Williams\\
  Department of Mathematics\\
  London School of Economics and Political Science\\
  London, UK\\
 \texttt{a.e.williams1@lse.ac.uk}}

\date{\empty}

\begin{document}
\maketitle

\begin{abstract}
With increases in population, there is a noticeable change across the world in pollution levels. Recently there has been growing demand for renewable energy operated devices boomed. Numerous reasons have led to such growth including lower operating costs and reduced greenhouse gas emissions. 

In order to obtain the optimised output, it is required to consider all the major parameters constraining the decision-making. Multi-Criteria Decision-Making (MCDM) is one of the most reliable and effective tools for decision making with many objectives. The technique focuses to prioritise the available alternatives in the decision space by considering the influential factors (or parameters) and their relative importance in the overall decision making. Because analysis conducted using MCDM approaches utilises an algorithm to obtain the desired output, this paper focuses on the application of an MCDM approach to identify those criteria that are essential in renewable energy technology systems. 

\vspace{2.0mm}
\end{abstract}

\section{Introduction}
Energy plays an important role in the economic development of a country, where the survey \cite{karanfil2015electricity} demonstrated a strong correlation between economic growth and electricity consummation. The overutilisation of fossil fuels has over time led to a significant negative impact on the environment. The report of the Intergovernmental Panel on Climate Change (IPCC) \cite{pachauri2014climate} states that by the end of the 21st century there may be an increase of 6 to 7 degrees in the earth’s temperature due to increased utilisation of fuels. Choosing alternate sources of energy is a multidimensional decision process that involves different parameters at different levels. 

The rises in both population and pollution levels means that there has been increased need for alternate sources of energy. In particular, global energy demand is expected to expand by around 1\% compoundly per annum until 20240 \cite{nejat2015global}. Further, present reports suggest that approximately one billion people across the globe do not have access to electricity \cite{gomez2020mathematical}. There has however been accelerated progress over the last few decades in moving to renewable sources of energy (Figure 1). Saini et al. \cite{saini2019design} demonstrate progress of PV wind microgrids that are specifically for electricity generation. Despite this, the power generated from renewable sources does not yet keep pace with the ever increasing growing demand. 

\begin{figure}[ht!]
\centering
	\includegraphics[width=0.85\linewidth]{{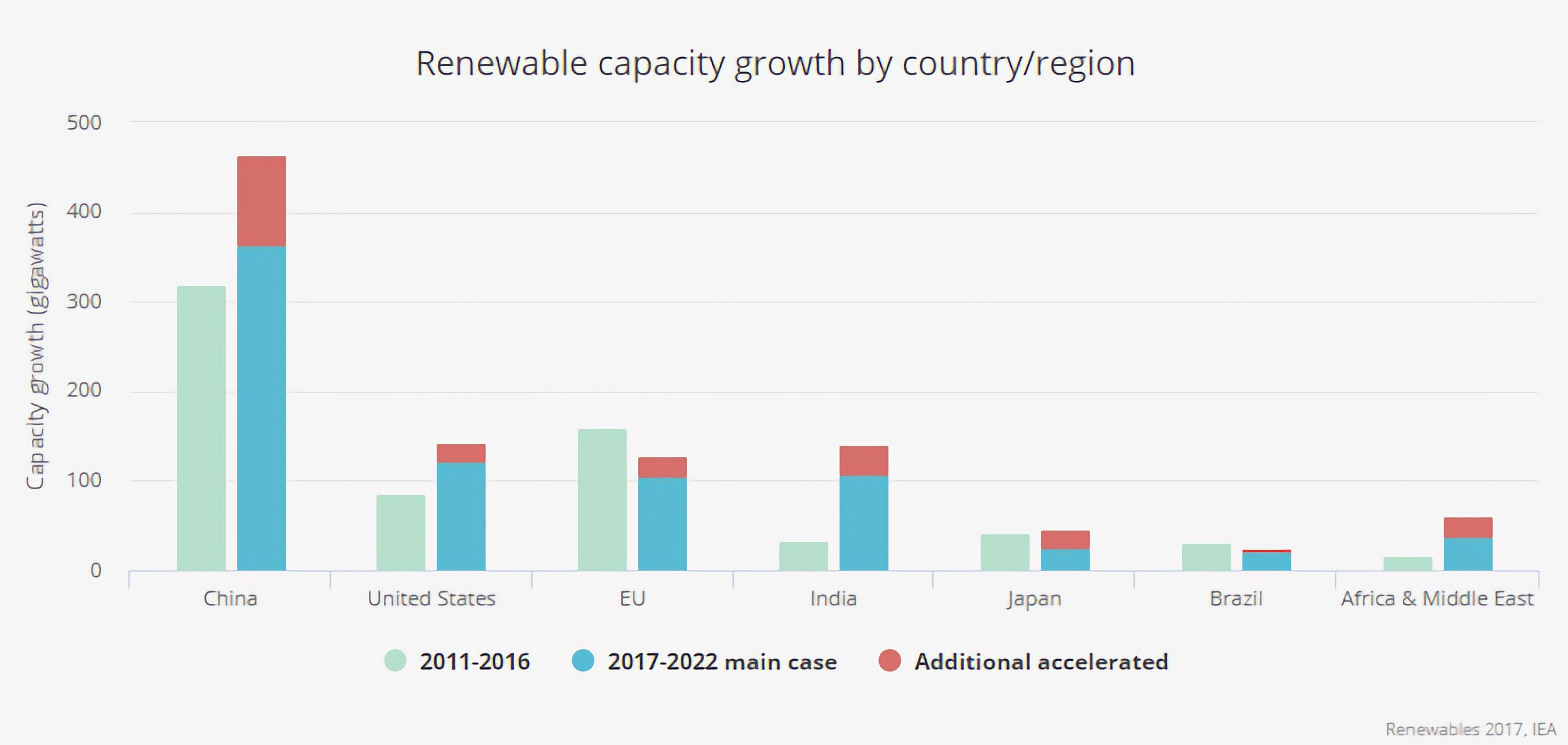}}
\caption{Adoption of renewable source of energy by different countries \cite{smith2018three}.}
    \label{fig:Figure1}
\end{figure}

In order to cope with the climate crisis of environmental sustainability and exhaustion of non-renewable energy, many countries are actively looking for renewable sources of energy production. This transition from traditional fossils to more clean sources of energy is a major challenge faced across the globe. 

This shift unsurprisingly impacts heavily on those countries whose economy depend upon the export of fossil fuels. Alizadeh et al. \cite{alizadeh2020improving}  provide a framework for a renewable energy system that includes the identification of the constraints which hinder the use of renewable sources of energy.

Upon designing a policy for renewable sources of energy, their low cost, sustainability and efficiency are key crucial considerations. Solar energy has proven to be one of the most efficient sources of electricity generation. The hazardous effect of conventional sources of energy on the environment has led to the development of green technology.  In particular, the development has focuses on looking into energy alternatives to reduce the emission of greenhouse gases \cite{ahmadi2018solar}.

Different renewable technologies are being explored in order to facilitate the increasing demand for sustainable energy \cite{guangqian2018hybrid}. The transport sector for example nearly emits of quarter of total carbon emissions \cite{Cozzi2023}), however, electric vehicles (EVs) have zero emission potential and some EV makers have seen year-on-year demand increases of around 185\%. 

With such surges in demand, it is necessary that much attention should be given to how such a global surge in demand for electricity is met. Bohanec et al. \cite{bohanec2017qualitative} provide a number of factors that are appropriate in order to select the energy production system.

An increase in the complexity and multiplicity involved in energy planning cannot be resolved by classical single objective optimisation techniques. In consequence, this places a major constraint on decision makers when optimising the available energy alternatives. Further, energy planning at set-up becomes more complicated due to topographical and geographical constraints \cite{tsoutsos2009sustainable}. Hence, MCDM tools were considered as an unbiased evaluation tool to resolve these constraints.  

Qian et al. \cite{qian2021fuzzy} implemented Fuzzy AHP and TOPSIS for the production location of Photovoltaic Energy and Solar Thermal Energy. Wang et al. \cite{wang2009review} performed an analysis based on the investment cost, C$O_2$ emissions efficiency, maintenance cost, land use and job creation. The purpose of the paper is to reduce the gap existing in literature by providing a set of parameters to develop a framework for the selection of solar panel by assigning a weight to evaluate the factors for standard dominance.

\section{Why Solar Energy?}
The ability of photovoltaic cells to convert the solar radiations directly into electricity make photovoltaics a future of the world of renewable energy technology. As the photovoltaics acquires the energy in decentralised way, it contributes in sustainable production of energy. In present scenario there are two types of PV panels i.e., monocrystalline and poly crystalline PV panels available for installation. Due to presence of more crystals embedded in monocrystalline panel, it is more efficient and quite expansive too. The statistics of 2022 indicates that (Figure 2) China is at the top in the production of renewable source of energy, using wind as one of its major sources. Europe achieve a first position in the production of solar energy. Data indicates that by as early as 2023, Europe could replace non-conventional source of energy by conventional energy.

\begin{figure}[ht!]
\centering
	\includegraphics[width=0.6\linewidth]{{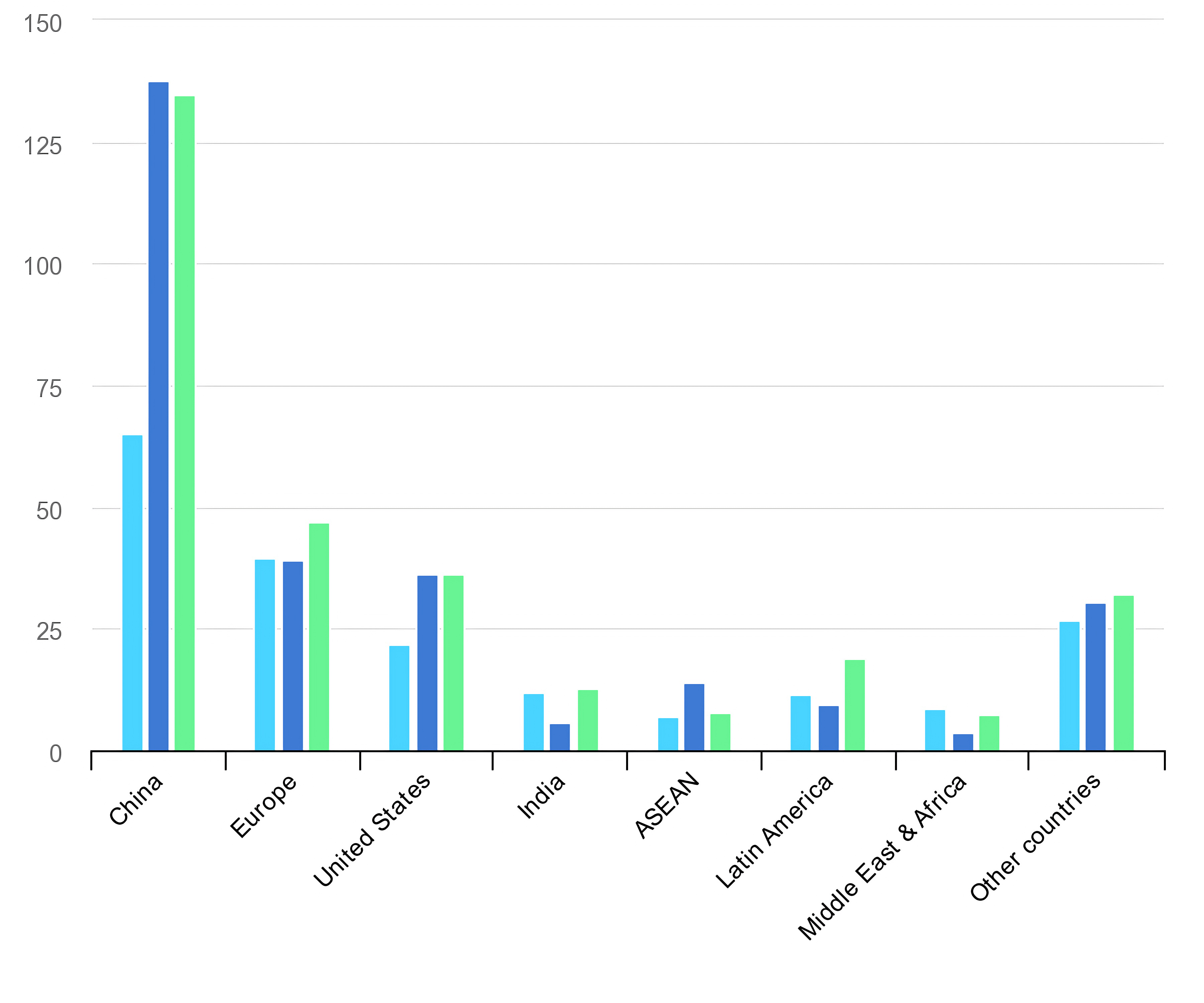}}
\caption{Energy Production across the Globe in 2022 \cite{douglasbroom2022}.}
    \label{fig:Figure2}
\end{figure}

The ever-increasing growing population and increasing environmental pollution made a paradigm shift of traditional energy sources to the renewable one (Figure 3). The installation of panels not only reduces the electricity cost but contributes to sustainable environment. Moreover, in present era people are switching to the electric vehicles due to the soar in the prices of fuels. Thus, the paper contributes significantly to identify the parameters of PV panels available in the market, their efficiency and develop a mathematical model to identify the optimality of these parameters. 

\begin{figure}[ht!]
\centering
	\includegraphics[width=0.6\linewidth]{{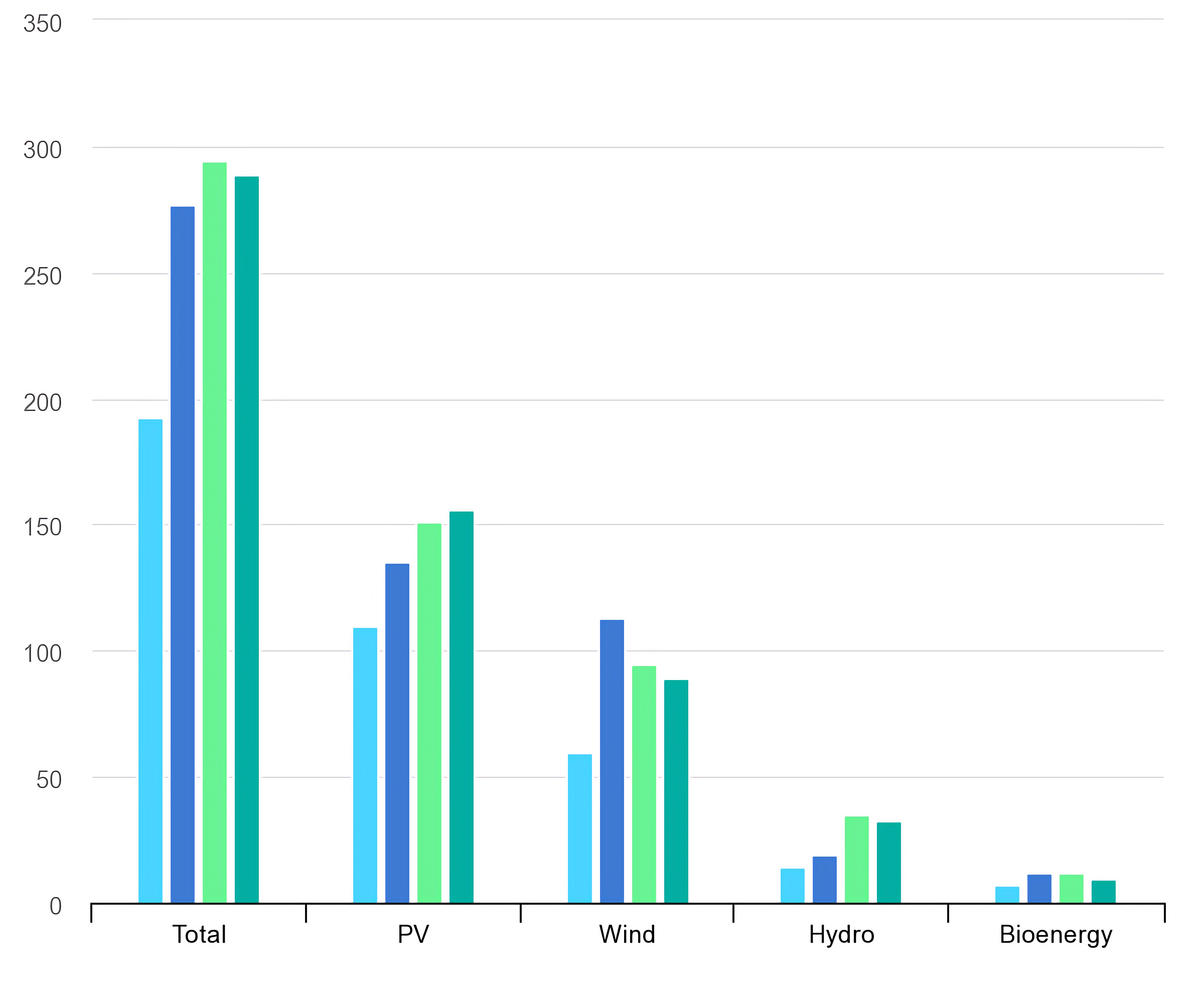}}
\caption{Renewable source of energy production in 2022 \cite{douglasbroom2022}.}
    \label{fig:Figure3}
\end{figure}

\newpage

\section{Mathematical Analysis}
Multi-Criteria Decision-Making (MCDM) is an appropriate tool in order to analyse the different criteria and alternatives defined in decision space. This is a mathematical approach that considers multiple criteria and provides a way to order them depending on their role in the decision making. MCDM approaches additionally provide tools to deal with multiple (perhaps) conflicting objectives, different forms of data and a large number of criteria. In particular, the methodology presented here provides a solution to increase the complex energy management system. 

Traditional programming approaches such as linear and integer linear programming identify the optimal output at minimum cost. In light of the inherent complexity of sustainability problems, there is a need to consider multiple objectives within one overarching framework. MCDM techniques are in consequence used to consider available criteria and alternatives defined in decision space.

\begin{figure}[ht!]
\centering
	\includegraphics[width=0.6\linewidth]{{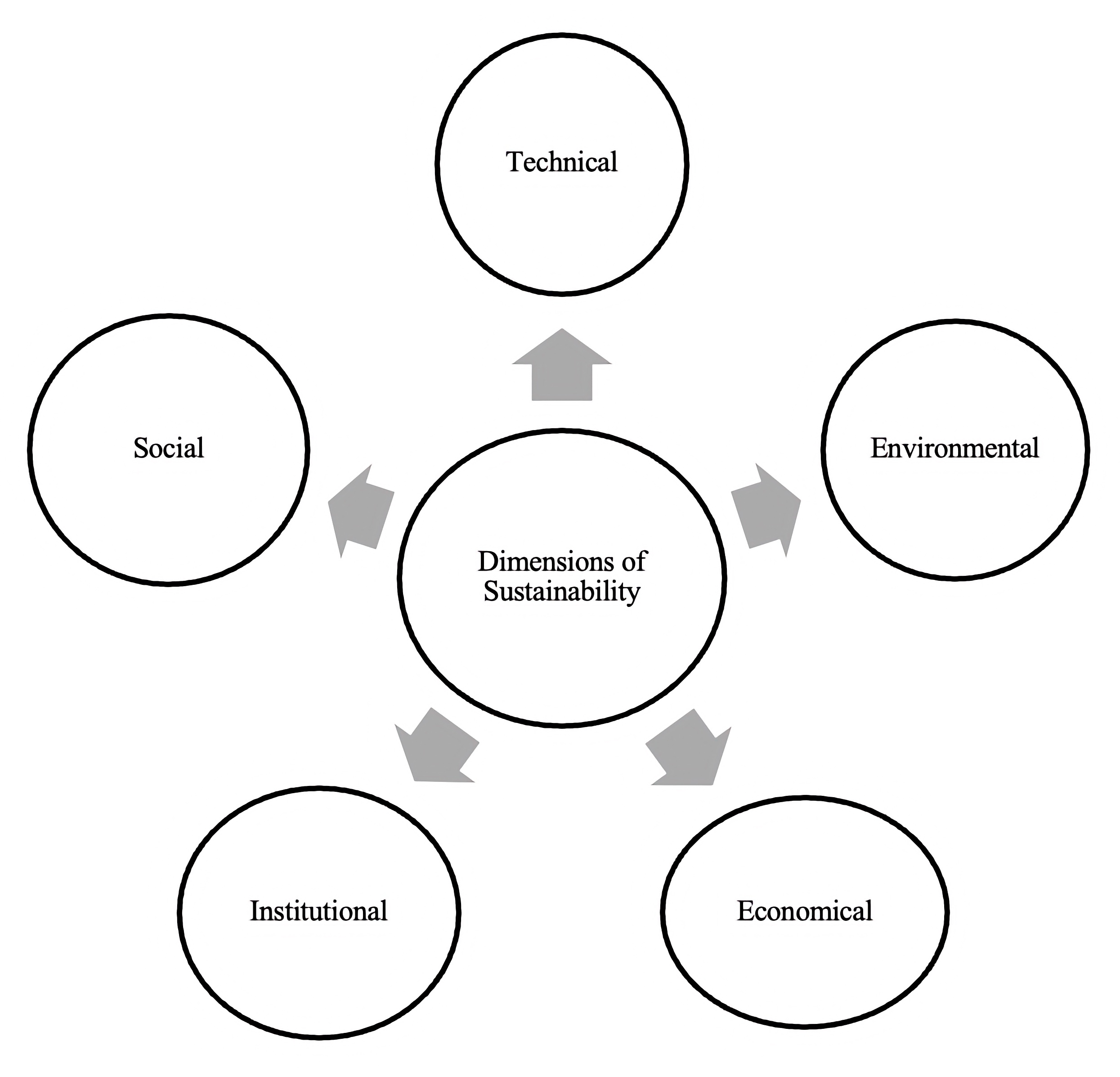}}
\caption{Key parameters of sustainability for electrification.}
    \label{fig:Figure4}
\end{figure}

Ilskog \cite{ilskog2008indicators} define five dimensions of sustainability, illustrated in Figure 4. The variance from the above sustainability parameters yields the inclusion of one additional parameter, namely a technical aspect that focuses to facilitate the proper electrification of the system. These indicators of sustainability are taken under consideration by reviewing the previous literature. 

This sustainability is then further ranked based on the sub-parameters \cite{wang2009review}, given below (Table 1). Many researchers have used these criteria in order to evaluate the sustainability of renewable systems.

\begin{table}[ht!]
\centering
\begin{tabular}{c|c|c|c|c|c|c}
Technical            &  & Economic         &  & Environmental        &  & Social          \\ \hline
Energy efficiency    &  & Investment cost  &  & Emission of gases    &  & Social benefits \\
Safety               &  & Operational cost &  & Emission of articles &  & Employability   \\
Reliability          &  & Fuel cost        &  & Land-use             &  &                 \\
Primary energy ratio &  & Electric cost    &  & Noise pollution      &  &                
\end{tabular}
\caption*{Table 1: Evaluation criteria for renewable source system}
\end{table}

Depending on the number of alternatives the differences can be highlighted between Multi Attribute Decision making (MADM) and Multi Objective Decision Making (MODM). MODM is an appropriate approach to evaluate continuous alternatives for which we define constraints in the form of vector decision variables. An optimal output is obtained by considering the set of alternatives that undergo the analysis by degrading the performance of one or more objectives. MADM makes \cite{azar2000multiattribute} the decision over the available alternatives that are characterized by multiple attributes defined in decision space. Different MCDM techniques are employed by researchers in the field of renewable energy to know the preference of the criteria and the available alternatives over the other. In this paper, the comparative analysis of different MCDM techniques is carried out in order to evaluate the optimal alternatives among the set of defined alternatives. 

More formally, MCDM techniques support the rational decision making which results in getting an optimized output by undertaking the problem constraints, namely
$$
\begin{aligned}
\text{minimise} \quad F(\boldsymbol{x}) = \left( f_1 (\boldsymbol{x}), f_2 (\boldsymbol{x}), \ldots, f_k (\boldsymbol{x}) \right)^T &\phantom{\ge} \\
\text{subject to } \quad g_i(\boldsymbol{x}) &\ge 0 \quad \text{for } i = 1,2, \ldots, n \\
x_i &\ge 0 \quad \text{for } i = 1,2, \ldots, n, 
\end{aligned}
$$
where $\boldsymbol{x} = (x_1, x_2, \ldots, x_n)^T \in \mathbb{R}^n$ is $n$-dimensional vector of decision variables. The inequalities are constraining the model for each $i$, while, the objective function is a vector of the parameters for the selection of attributes of PV panels. 

Hence, MCDM tools have become a popular tool in energy planning since it supports the decision-makers to formally evaluating criteria defined in the decision space and rank them based on the decided priority in the space. In order to maximise the probability that a correct decision is made and to provide an optimal energy plan it is preferable to take a large set of criteria that are defined in some decision space. Intuitively, the larger the set of criteria used will increase the likelihood that a correct decision is made, however, this will be at the cost of  significant increases in computation time. Further, MCDM techniques can help the decision maker quantify the criteria depending on its importance in decision making.  

This research work gives a brief about different MCDM techniques and provides insight into energy planning based on renewable energy sources, where we focus primarily on solar energy. The MCDM techniques formulated in the research work used to find out the optimal solution to the energy system design involving multiple and (perhaps) conflicting objectives. Ranking of technologies in order to find the most optimized factor that impact the decision-making influence the relative importance of the determinant importance of the parameters in decision making.

\section{Selection Criteria}
The standard deviation and entropy method (Table 2) is used for assigning a weight to different criteria defined in decision space. The weights are assigned to some criteria before then evaluating their importance in the decision making. The value of the weights assigned to the criteria depends on the nature of the given criterion, the dispersion of choice performance, the objective uncertainty and the analysis done by decision makers when deal with small differences in alternatives. 

The standard deviation method determines the target related to the alternatives which are defined in decision space \cite{diakoulaki1995determing}. Depending on the benefit and cost targets, a normalised matrix is then formulated. For each objective, the standard deviation is calculated \cite{keshavarz2021determination}.

\begin{table}[ht!]
\centering
\begin{tabular}{c|c|c}
Attributes                 & Standard Deviation Method & Entropy Method \\ \hline
PV efficiency              & 0.040540                  & 0.077422       \\
PV lifetime                & 0.0767126                 & 0.011123       \\
Total power generation     & 0.023082                  & 0.027539       \\
PV panel cost              & 0.009138                  & 0.487074       \\
Battery cost               & 0.0758723                 & 0.384587       \\
Hourly self-discharge rate & 0.774653                  & 0.012256      
\end{tabular}
\caption*{Table 2: Weights assigned to select parameters using the standard deviation and entropy methods.}
\end{table}

Figure 5 below illustrates that for the selection of PV panel all the parameters \enquote{significantly} contribute to decision making. Weights is assigned to the attributes to remove the biasness in decision making and to get an optimal output. The graph indicates that there are significant differences in few parameters in the selection of solar panels. Note that some of the parameters such as generation of power does not show the deviation in the assigned weights. 

Figure 5 indicates the weights assigned to the decision parameters chosen to select the PV panel depending on the different attributes. The weights assigned to the attributes evaluate the data on the same scale resulting in the reduction of bias.  

\begin{figure}[ht!]
\centering
	\includegraphics[width=0.75\linewidth]{{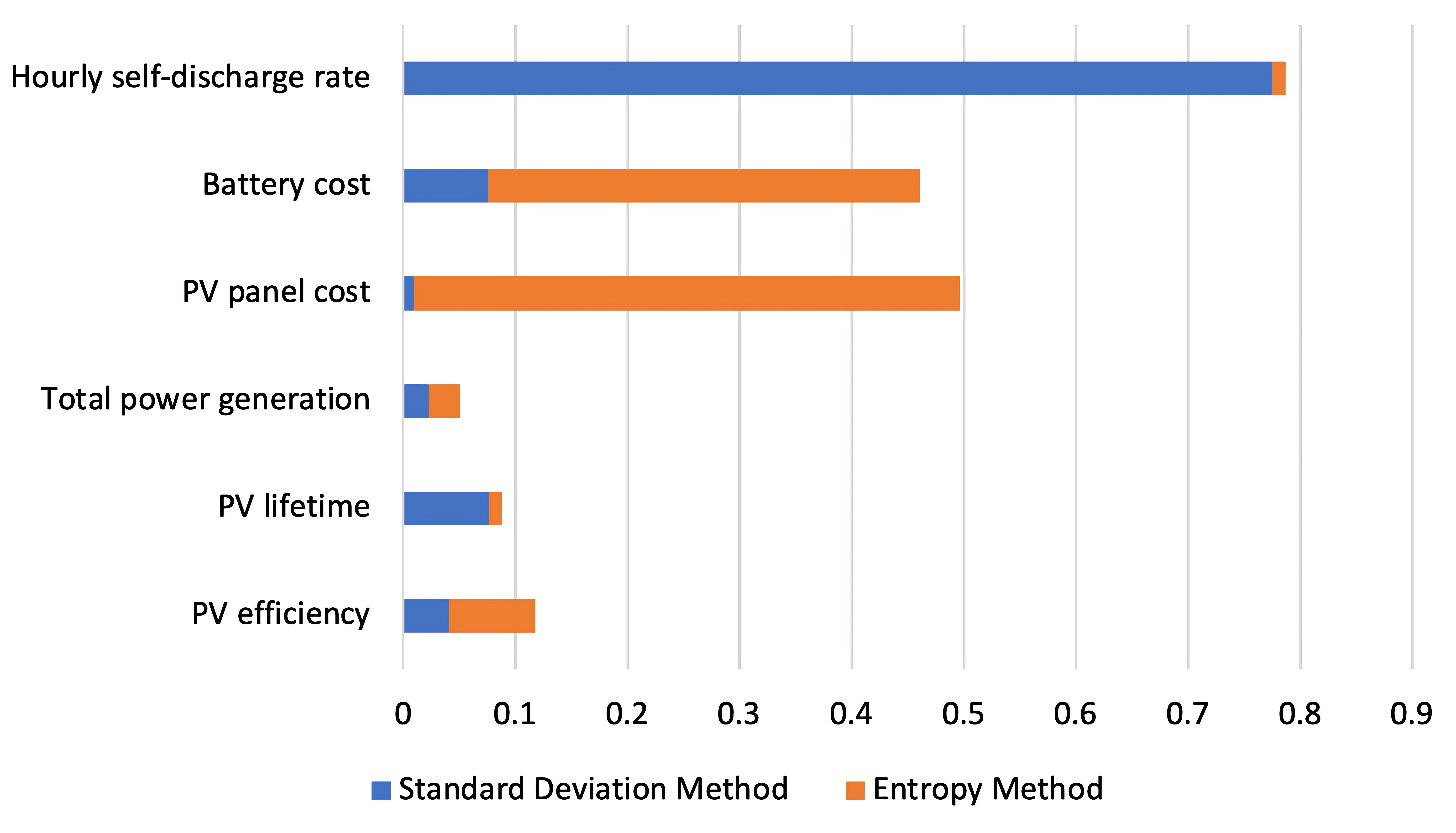}}
\caption{Weights assigned to attributes.}
    \label{fig:Figure5}
\end{figure}

In this paper TOPSIS and MOORA are used to check the accessibility of the results. TOPSIS is one of the MCDM tool that rank the alternatives by calculating the distance from both positive and negative ideal solution simultaneously. The optimal alternative is determined by the highest closeness coefficient. According to TOPSIS algorithm the normalized matrix is formulated. 

\begin{table}[ht!]
\centering
\begin{tabular}{ccccc|cc}
\multicolumn{5}{c}{TOPSIS}                           & \multicolumn{2}{|c}{MOORA} \\ \hline
Alternatives & Si+      & Si-      & Ci       & Rank & Augmented Matrix  & Rank  \\ \hline
A1           & 0.022217 & 0.020747 & 0.482891 & 25   & -0.1428           & 23    \\
A2           & 0.012274 & 0.022345 & 0.645451 & 7    & -0.13549          & 8     \\
A3           & 0.012967 & 0.019885 & 0.605284 & 11   & -0.13812          & 15    \\
A4           & 0.008079 & 0.024527 & 0.752216 & 3    & -0.13026          & 4     \\
A5           & 0.012585 & 0.019926 & 0.612896 & 10   & -0.13661          & 13    \\
A6           & 0.022283 & 0.016468 & 0.424979 & 29   & -0.15055          & 29    \\
A7           & 0.013431 & 0.021769 & 0.618429 & 8    & -0.13642          & 11    \\
A8           & 0.01418  & 0.017799 & 0.556584 & 16   & -0.14056          & 18    \\
A9           & 0.016049 & 0.019799 & 0.552307 & 18   & -0.14046          & 17    \\
A10          & 0.01729  & 0.017507 & 0.503112 & 23   & -0.14313          & 25    \\
A11          & 0.016688 & 0.019104 & 0.533753 & 20   & -0.14189          & 22    \\
A12          & 0.006458 & 0.024918 & 0.79417  & 2    & -0.12823          & 2     \\
A13          & 0.017691 & 0.013588 & 0.434419 & 27   & -0.14569          & 27    \\
A14          & 0.015075 & 0.017349 & 0.535054 & 19   & -0.14131          & 20    \\
A15          & 0.01603  & 0.015786 & 0.496174 & 24   & -0.14292          & 24    \\
A16          & 0.013675 & 0.019851 & 0.5921   & 14   & -0.13784          & 14    \\
A17          & 0.011343 & 0.02111  & 0.650492 & 6    & -0.13516          & 6     \\
A18          & 0.012627 & 0.020347 & 0.617062 & 9    & -0.13605          & 9     \\
A19          & 0.013601 & 0.020668 & 0.603106 & 12   & -0.1354           & 7     \\
A20          & 0.003054 & 0.027819 & 0.901089 & 1    & -0.12368          & 1     \\
A21          & 0.022062 & 0.012702 & 0.365369 & 30   & -0.15109          & 30    \\
A22          & 0.015799 & 0.019561 & 0.553203 & 17   & -0.13984          & 16    \\
A23          & 0.010791 & 0.023432 & 0.684676 & 5    & -0.13173          & 5     \\
A24          & 0.016133 & 0.023225 & 0.590093 & 15   & -0.13661          & 12    \\
A25          & 0.014552 & 0.021532 & 0.596723 & 13   & -0.13634          & 10    \\
A26          & 0.01991  & 0.018212 & 0.47773  & 26   & -0.14431          & 26    \\
A27          & 0.017602 & 0.018994 & 0.519016 & 21   & -0.14084          & 19    \\
A28          & 0.007796 & 0.023441 & 0.750422 & 4    & -0.12927          & 3     \\
A29          & 0.020914 & 0.0158   & 0.430357 & 28   & -0.14697          & 28    \\
A30          & 0.02023  & 0.02074  & 0.506227 & 22   & -0.14162          & 21   
\end{tabular}
\caption*{Table 3: Rank assigned to the alternatives by TOPSIS and MOORA.}
\end{table}

\section{Results and Discussions}

In this paper TOPSIS and MOORA are used to check the accessibility of the results. TOPSIS is one of the MCDM tool that rank the alternatives by calculating the distance from both positive and negative ideal solution simultaneously. The optimal alternative is determined by the highest closeness coefficient. According to TOPSIS algorithm the normalized matrix is formulated. 

\begin{figure}[ht!]
\centering
	\includegraphics[width=0.75\linewidth]{{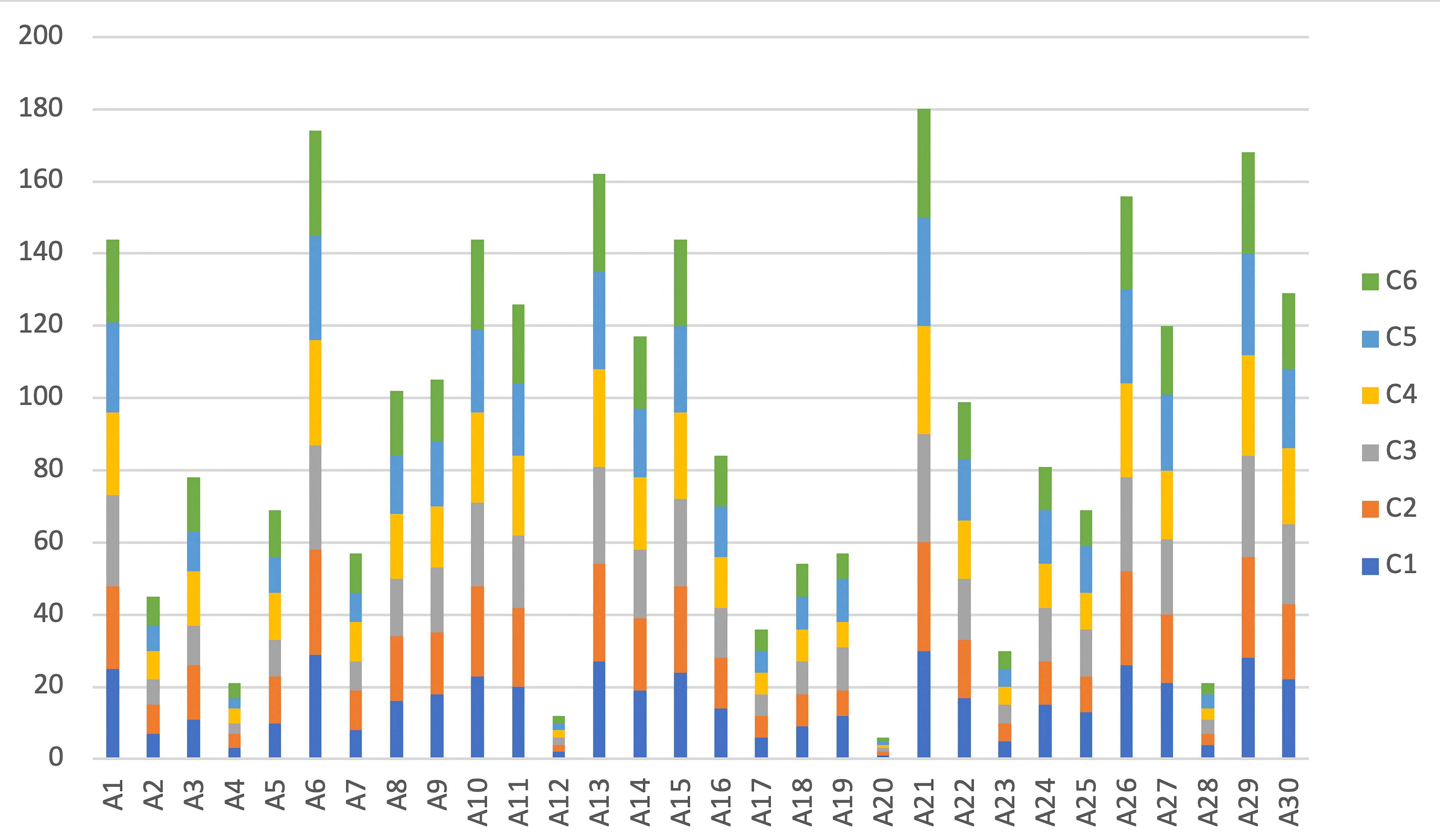}}
\caption{Comparison of MCDM Technique.}
    \label{fig:Figure6}
\end{figure}

The competency of conventional sources depends directly on their geographical location and the ambient weather conditions. Thus, for the evaluation of such parameters, the application of MCDM techniques is the most appropriate tool. These approaches are often facilitated to choose the correct alternatives which are conflicting in nature and available for the decision space depending on different criteria. In order to adopt the most optimized alternative inter and intra-comparison of techniques is being carried out \cite{lu2007multi}.

\begin{figure}[ht!]
\centering
	\includegraphics[width=0.75\linewidth]{{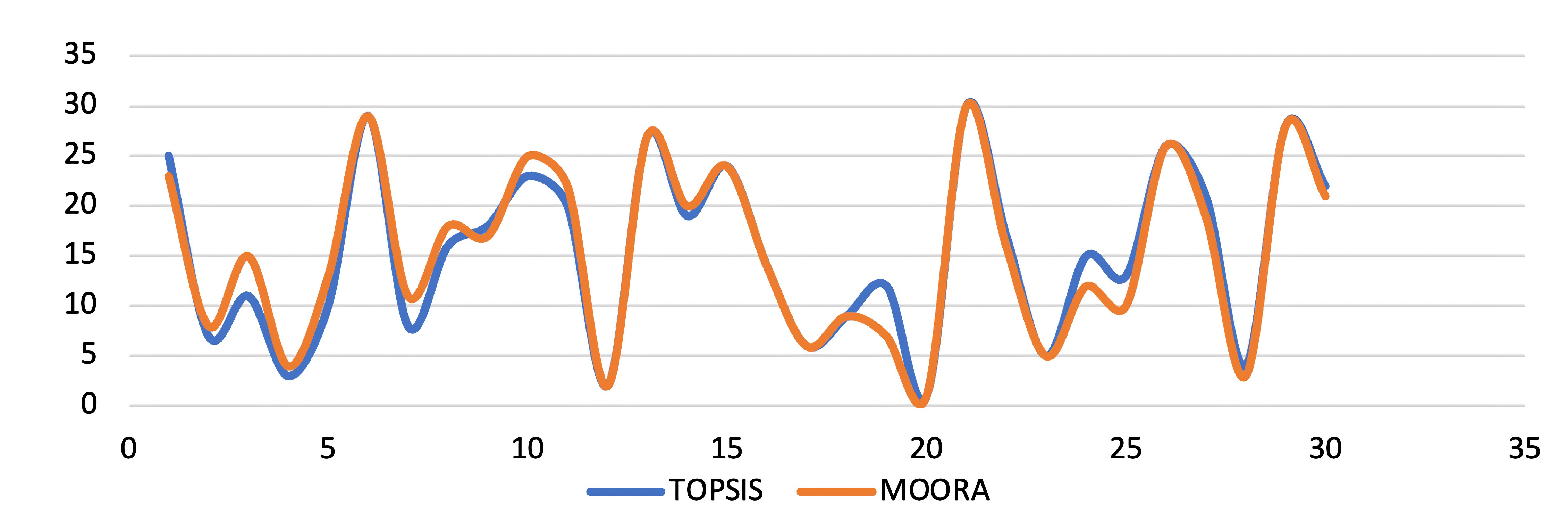}}
\caption{Comparison of MCDM Technique.}
    \label{fig:Figure7}
\end{figure}

In comparison with more conventional classical sources of energy, renewable energy provides a clean and theoretically inexhaustible form of energy. Despite this, there are obstacles that make it difficult for renewable energy sources to currently complete with conventional sources. It is for example more challenging to manage the variable supply from renewable sources in order to meet  demand. Further, for example the sunk costs associated with a renewable system can be higher in some geographical areas. 

For solar PV modelling, researchers have previously considered the constant parameters such as ideality factor, series and shunt resistance. The differences in the results were observed as all these parameters changes with temperature and insulation \cite{bana2016mathematical}. The result of the model indicates (Figure 7) that both the techniques rank the same alternative $A_{20}$ as 1 depending on the parameters defined in decision space.  It can be observed that the highest impact criteria impacting the preference value is $C_5$ (discharge rate of PV panel) and the second criterion in terms of local importance is lifetime of PV panel (Figure 5). Thus, it gives a valid explanation why the MCDM techniques order the alternative $A_{20}$ as 1. For entropy method the weights assign to different criteria as compared to that of standard deviation method. The criterion of PV panel cost and battery cost involved in decision making of optimal solution.


\newpage

\bibliographystyle{plain}
\bibliography{references}

\end{document}